\documentclass{amsart}

\usepackage{amssymb, latexsym}
\usepackage{graphicx}
\title{A Reciprocity Theorem for Monomer-Dimer Coverings}

\author{Nick Anzalone \and John Baldwin \and Ilya
Bronshtein \and T. Kyle Petersen}


\newtheorem{theorem}{Theorem}
\newtheorem{lemma}{Lemma}
\newtheorem{defin}{Definition}

\begin{document}

\maketitle

\begin{abstract}
The problem of counting monomer-dimer coverings of a finite patch in a
lattice is a longstanding problem in statistical mechanics.  It has only 
been exactly solved for the special case of dimer coverings in two
dimensions (\cite{kast}, \cite{tf}). In earlier work, Stanley
(\cite{stanley}) proved a reciprocity principle governing the number
$N(m,n)$ of dimer coverings (also known as perfect matchings)
of an $m$ by $n$ rectangular grid, where $m$ is fixed and $n$ is
allowed to vary. As reinterpreted by Propp (\cite{propp}), Stanley's
result concerns the unique way of extending $N(m,n)$ to $n < 0$ so
that the resulting bi-infinite sequence, $N(m,n)$ for $n \in
\mathbb{Z}$, satisfies a linear recurrence relation with constant
coefficients. In particular, Stanley shows that $N(m,n)$ is always
an integer satisfying the relation $N(m,-2-n) =
\epsilon_{m,n}N(m,n)$ where $\epsilon_{m,n} = 1$ unless $m\equiv$
2 (mod 4) and $n$ is odd, in which case $\epsilon_{m,n} = -1$.
Furthermore, Propp's method was applicable to higher-dimensional
cases, such as the dimer model on an $m_1$-by-$m_2$-by-$n$ grid.
This paper discusses similar investigations of the number
$M(m,n)$, of monomer-dimer coverings, or equivalently, (not
necessarily perfect) matchings of an $m$ by $n$ rectangular grid.
We show that for each fixed $m$ there is a unique way of extending
$M(m,n)$ to $n < 0$ so that the resulting bi-infinite sequence,
$M(m,n)$ for $n \in \mathbb{Z}$, satisfies a linear recurrence
relation with constant coefficients.  We show that $M(m,n)$,
\emph{a priori} a rational number, is always an integer, using a
generalization of the combinatorial model offered by Propp.
Lastly, we give a new statement of reciprocity in terms of
multivariate generating functions from which Stanley's result
follows.
\end{abstract}

\section{Introduction}
\label{sec:in}

\subsection{Background}
The problem of counting the monomer-dimer coverings of a finite
patch of a lattice has been examined for many years in the field of
statistical mechanics (see \cite{krs}), and has applications in
biology, chemistry and physics.  The closely related dimer
problem, where the number of monomers is zero, was exactly solved
for two dimensional lattices by Kasteleyn (\cite{kast}) and
Temperley and Fisher (\cite{tf}).  While we cannot offer an exact
solution to the monomer-dimer problem, even in two dimensions, we can
describe a symmetry property that  must be satisfied by exact
solutions to various restrictions of this problem, in which one
considers a sequence of subgraphs of a lattice, where the successive
subgraphs grow in one direction  while their size in all other
directions remains fixed.

In particular, we consider families of subgraphs of lattices indexed 
by the natural numbers ($n = 1,2,3,\ldots$), and count the number
$a_{n}$ of monomer-dimer coverings of the $n$th subgraph.
The numbers $a_{n}$ satisfy a linear recurrence with constant 
coefficients
and therefore are given by an exact formula of polynomial-exponential 
type; 
one can thus in each individual case write down an exact formula for 
$a_n$, 
at least in principle, and in particular, one can study symmetries of 
the 
resulting function of $n$ when $n$ is no longer restricted to being a 
natural number.  Alternatively, when $n$ is a negative integer,
one can often find $a_n$ by working directly with the recurrence
relation satisfied by the sequence.  Surprisingly, we find
that the new numbers are related, term by term, to the original
sequence of numbers; for instance, in some cases $a_{n} \sim a_{-n-2}$.  
Richard Stanley's book \cite{stan1}, 
in the context of rational generating functions, devotes 
an entire section to exploring the relationships (called 
\emph{reciprocity}
relationships) between positively- and nonpositively-indexed terms
of a sequence.  The nonpositively-indexed terms may even have a nice
combinatorial meaning on their own, as seen in the case of Ehrhart
reciprocity, for example.\footnote{Ehrhart reciprocity describes a 
relationship between the number of lattice points found in a closed 
rational polytope, and the number of lattice points found in its 
interior.  
See \cite{stan1}.} 
In this paper we describe a relationship between terms of certain 
integer
sequences $\{a_{n}\}$, $n \in \mathbb{Z}$, as well as give a 
combinatorial
interpretation to the ``un-natural'' terms.

\subsection{The Problem}

Rather than use the language of statistical mechanics (monomers
and dimers on lattices), we will refer to matchings of graphs
where we define a \emph{matching} of a graph $G = (V,E)$ to be a
collection of edges of $G$, no two of which share a vertex, together
with all the vertices of $G$ that are not incident with those
edges.\footnote{Another equivalent terminology would be to speak
of tilings of planar regions where the tile set consists of a
domino (2 by 1 rectangle) and a monomino (1 by 1 rectangle).}  Edges
correspond to dimers, isolated vertices correspond to monomers.  A
matching is called \emph{perfect} if it is composed entirely of
edges, as in the dimer problem. If we consider the number of
matchings of a rectangular grid-graph of fixed height $m$, and
varying width $n$, we obtain an integer sequence, $\{a_{n}\}$. For
example, when $m = 1$, we obtain the Fibonacci sequence:
\[ 1,2,3,5,8,\ldots \]  This sequence satisfies a linear
recurrence, namely $a_{n} = a_{n-1} + a_{n-2}$. So we can always
obtain the value of a term based on the two terms to the left of
it. But similarly, we can obtain the value of a term by the two
terms to the right of it, e.g. $3 = 8-5$, or $a_{n} = a_{n+2} -
a_{n+1}$.  In this way we can extend the Fibonacci sequence to the
left and get values for $a_{n}$ when $n \leq 0$:
\[ \ldots -8,5,-3,2,-1,1,0,1,1,2,3,5,8, \ldots \]  
So we see that we now have a doubly infinite integer sequence, and that 
it
is symmetric up to sign, i.e. $a_{n} = \pm a_{-n-2}$.  
Some natural questions arise: Does this symmetry mean something?  
Do the values of $a_{n}$ for $n < 0$, being integers, count something
(at least up to sign)?
If so, can we extend the result to larger values of $m$?

The answer to all the questions is yes.  In \cite{propp}, Propp
considered integer sequences $N(m,n)$, generated by perfect
matchings of rectangular grid graphs of fixed height, and how they
extend to values for $n < 0$.\footnote{His first nontrivial case,
$N(2,n)$ is also the Fibonacci sequence.}  He came up with a
unified combinatorial object and a way of counting signed
matchings of these objects that allowed him to show that any
half-infinite sequence given by the number of perfect matchings of
a fixed-height grid graph extends to a bi-infinite sequence that
does three things:
\begin{itemize}
\item It satisfies a linear recurrence relation of finite degree 
with constant coefficients.
\item The numbers obtained when going backwards are unique and are
always integers.
\item The bi-infinite sequences have a special
kind of symmetry, or \emph{reciprocity}, stated roughly as
$|N(m,n)| = |N(m,-n-2)|$ for any $m > 0$.
\end{itemize}

Matchings in general (of the graphs we will classify) satisfy the
same sort of linear recurrence (as we will explain in section 3).
But even given that a covering problem satisfies a linear
recurrence, integrality is not ensured when running the recurrence
backwards.\footnote{\emph{A priori} the numbers obtained from
running a linear recurrence in reverse need only be rational.  If
the recurrence is of the form $a_{n} = c_{1}a_{n-1} + c_{2}a_{n-2}
+ \cdots + c_{k}a_{n-k}$, then to push the sequence backwards we
solve for the term of smallest index to get $a_{n-k} = \frac{a_{n} -
c_{1}a_{n-1} - \cdots - c_{k-1}a_{n-k+1}}{c_{k}}$.} For a simple
example, consider covering a 2-by-$n$ grid with monomers and
dimers, but only allow the dimers to be vertical.  (It is easy
to see that the number of coverings is $2^n$, but in keeping
with the spirit of this article, we ignore the exact formula
and work with a linear recurrence instead.)  The number of
such coverings, $a_n$, is governed by $a_{n} = 2a_{n-1}$.  Upon 
observing that $a_{1} = 2$, we can generate the half-infinite sequence
$2,4,8,16,\ldots$.  To work the sequence backwards, we reverse the
recurrence by writing it as $a_{n} = \frac{1}{2}a_{n+1}$ and then 
generate the bi-infinite sequence \[\ldots \frac{1}{8}, \frac{1}{4},
\frac{1}{2}, 1, 2, 4, 8, \ldots\] For $n = -1$ we cease to observe
integrality.
However in the case of monomer-dimer tilings, we can guarantee
integrality of the sequences by showing they can be generated by
counting matchings of certain graphs for \emph{all} values of $n$.

The objects we will present not only ensure integrality, but also give 
a
direct proof of the paper's main idea: a reciprocity statement for
general matchings.  This claim would appear to be the least obvious
judging by the numerical evidence. Indeed, consider the integer
sequence $M(2,n)$, the number of matchings of a 2-by-$n$
rectangular grid graph. The sequence is governed by the recurrence
\[M(2,n) = 3M(2,n-1) + M(2,n-2) - M(2,n-3)\] with initial
`natural' values $M(2,1) = 2$, $M(2,2) = 7$, and $M(2,3) = 22$,
allowing us to generate the bi-infinite sequence
\[ \ldots 14,11,2,3,0,1,0,1,2,7,22,71,228, \ldots \] Although we get
integers when going backwards, we seem to lack the symmetry that
Stanley and Propp observed for perfect matchings.  But a premonition of
symmetry can be seen if one considers the parity (odd versus even) of each
term; one finds that $M(2,n)$ and
$M(2,-2-n)$ are congruent mod 2.

Why there should be symmetry (or symmetry mod 2) in integer sequences 
as above is explained by the answer to the question: what do the 
numbers,
$M(m,n)$ for $n<0$, count?  Propp had an answer for the case of
perfect matchings, but it does not directly translate to general
matchings.  Section 2 outlines the dilemma and gives a solution by
way of objects called ``empty vertices". In 
\cite{propp}, weights were not given to isolated vertices, simply
because there were none.  But when examining matchings as we have
defined them (a monomer-dimer model), isolated vertices do have
weights, giving rise to the need for empty vertices to explain the
observed phenomenon and tell the ``proper'' combinatorial story.
By incorporating empty vertices we will obtain a combinatorial
model whose (signed) number of \emph{signed} matchings is $M(m,n)$ 
for all $n$. Looking at the signed matchings we will be able to see 
the symmetry hidden in the terms above: \[\ldots 121-107, 41-30, 
12-10, 5-2, 1-1, 1, 0, 1, 1+1, 5+2, 12+10, 41+30, 121+107, \ldots\]
This also explains the symmetry mod 2.

The properties exhibited by the numbers $M(m,n)$ and the objects
they count can also be seen in terms of the generating function
$F_{m}(t,x,y,z) = \sum_{n=1}^{\infty}f_{n}(x,y,z)t^{n}$.  Here $x$
is the weight given to horizontal edges, $y$ is the weight given
to vertical edges, and $z$ is the weight given to vertices. The
polynomial $f_{n}(x,y,z)$ encodes all of the matchings of an $m$
by $n$ grid graph. We will show
\begin{enumerate}
\item $F_{m}(t,1,1,1) = \sum_{n=1}^{\infty}M(m,n)t^{n}$ \item
$F_{m}(t,1,1,0) = \sum_{n=1}^{\infty}N(m,n)t^{n}$ \item
$F_{m}(t,x,y,z) \sim -F_{m}(1/t,x,-y,-z)$
\end{enumerate}
The second item in the list shows how perfect matchings may be
extracted from general matchings. The last item in the list is a
heuristic expression of reciprocity that will be made precise in
section 4.2.

The theory of empty vertices is applicable to more
general situations than simply matchings of rectangular grid
graphs.  We will prove results for matchings of graphs of the form
$G \times P_{n}$, where $G$ is an arbitrary finite graph and
$P_{n}$ is the path graph of length $n$.  This is the ``box-product''
of the two graphs: a vertex of $G \times P_n$ is a pair $(u,v)$
consisting of a vertex of $G$ and a vertex of $P_n$, and two
such vertices $(u,v)$, $(u',v')$ are connected by an edge iff
either $u=u'$ and $v,v'$ are adjacent in $P_n$ or $v=v'$ and
$u,u'$ are adjacent in $G$.  We will refer to such
graphs as ``generalized rectangles". Notice that an $m$ by $n$
rectangular grid graph is $G\times P_{n}$ where $G = P_{m}$. If
$M(n)$ is the number of matchings of $G \times P_{n}$, we can 
show that there are objects (signed graphs) for which the (signed)
number of matchings is $M(-n)$, i.e., we can give a definition 
for $G\times P_{-n}$. Reciprocity asserts that the (signed) number of signed matchings
of $G \times P_{-2-n}$ is equal to (or equal to the negative of) 
the (unsigned) number of matchings of $G \times P_n$.
Throughout most of the paper, our proofs will be for
generalized rectangles, but examples will usually involve grid
graphs.

\section{Signed Graphs and Signed Matchings}
\label{sec:sign}

\begin{figure} [h]
\centering
\includegraphics{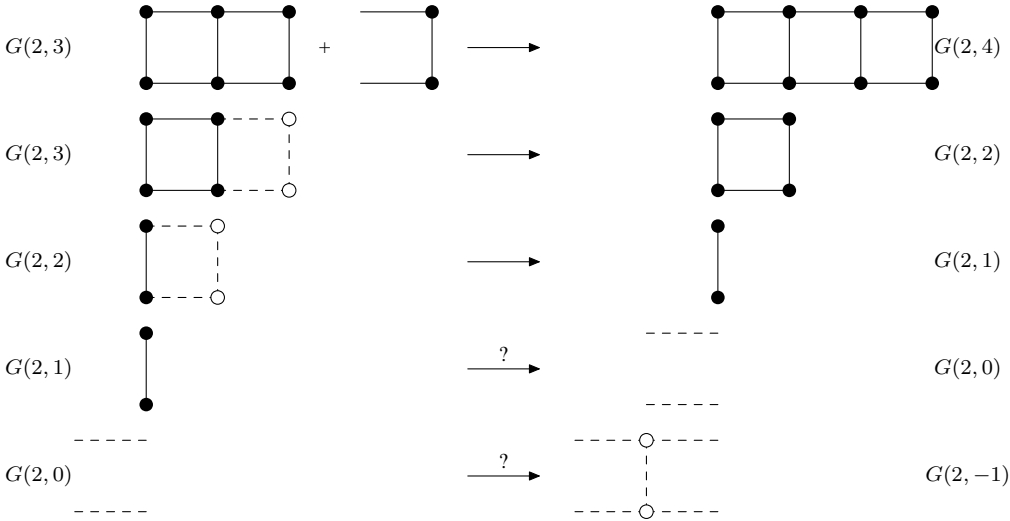}
\caption{Moving from $G(2,3)$ to $G(2,4)$ or $G(2,2)$}
\label{fig:move}
\end{figure}

Let us try to deduce the nature of $m$ by $n$ rectangular grid
graphs for all $n$, not concerning ourselves with matchings at
all.  In Figure \ref{fig:move}, denote the 2 by 3 grid graph by
$G(2,3)$.  The graph $G(2,4)$ is obtained from $G(2,3)$ by adding
two horizontal edges, a vertical edge, and two vertices as shown.
Likewise, we can obtain $G(2,2)$  by
removing the same set of edges and vertices from the right of
$G(2,3)$.  Then we get $G(2,1)$ and $G(2,0)$ inductively. But what
are we to make of $G(2,0)$?  It has no vertices, and two ``anti-"
horizontal edges. How many matchings should this have? What about
$G(2,-1)$? $G(2,-2)$?

\begin{figure} [h]
\centering
\includegraphics{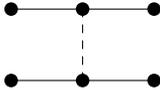}
\caption{Propp's model for $G(2,-1)$.}
\label{fig:propp}
\end{figure}

The definitions that Propp created for $G(m,n)$ with $n<0$ look
quite similar to the pictures in Figure \ref{fig:move}, except
that his horizontal lines are solid and there are vertices on the
endpoints of the horizontal edges where there are none in Figure
\ref{fig:move} (see Figure \ref{fig:propp}).

\begin{figure} [h]
\centering
\includegraphics{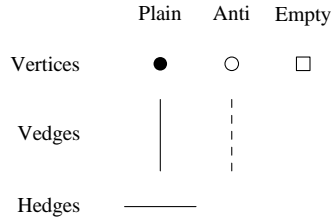}
\caption{Basic components of a signed graph.} \label{fig:comp}
\end{figure}

Notice that in Figure \ref{fig:move} the graphs for $G(2,0)$ and
$G(2,-1)$ have some horizontal edges that are incident with only
one vertex. The meaning of these missing vertices, these places
that are ``empty" of vertices, is described below.  Let $G =
(V,E)$ be a graph in the usual sense, i.e. a set of vertices and
edges between them, except that there is more than one kind of
vertex and more than one kind of edge (see Figure \ref{fig:comp}).
There are plain vertices, anti-vertices, and empty vertices, 
as well as plain vertical edges (vedges), anti-vedges, and
plain horizontal edges (hedges).
For now, plain components are given weight 1, anti-vertices and
anti-vedges are given weight -1, and empty vertices have weight
zero.\footnote{Later we will attach formal variables to edges and
vertices, but empty vertices will still have weight zero. Thinking
of monomer-dimer coverings, empty vertices represent positions in
the lattice that may not be occupied by a monomer.}

Define a \emph{signed matching} of a graph $G$ to be a collection
of non-adjacent edges of $G$ and all the vertices not incident
with those edges.  The weight of the matching is the product of
the weights of the components in the matching.  Note that this
implies that a matching that involves an edge one or both of
whose endpoints is an empty vertex has vanishing weight.
Matchings of weight zero ``do not count'' for purposes
of weighted enumeration, so in practice
it helps to think of empty vertices as being \emph{needy}.  For
any matching of nonzero weight, empty vertices \emph{need} to 
have one of their incident edges included.

\begin{figure} [h]
\centering
\includegraphics{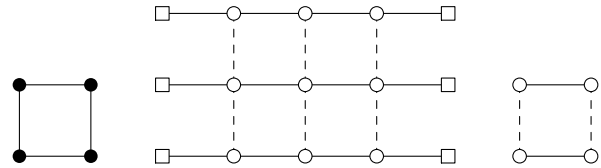}
\caption{Pictures of graphs $G(2,2)$, $G(3,-3)$, and
$G^{*}(2,2)$.} \label{fig:graph}
\end{figure}

Define a regular signed graph, $G(m,n) = P_{m} \times
P_{n}$ in the following way.  If $n > 0$, $G(m,n)$ is just a
rectangular grid graph, i.e., a graph with $mn$ plain vertices
arranged into $m$ rows of $n$ vertices each, with plain edges
adjoining horizontal and vertical neighbors.  Define the conjugate
graph of $G$, denoted $G^{*}$, to be the graph obtained by
replacing all of the vertices and vedges of $G$ with their
anti-counterparts (assuming that an anti-anti-vertex is a plain
vertex, etc.) but leaving the hedges alone. Under this definition, 
$(G(m,n))^{*} = G^{*}(m,n)$
is a graph with $mn$ anti-vertices arranged into $m$ rows of $n$
vertices with hedges adjoining horizontal neighbors and
anti-vedges adjoining vertical neighbors. 
Then for $n \leq 0$, $G(m,n)$ is defined to be a copy of $G^{*}(m,n)$ 
with a column of $m$ empty vertices on the left and $m$ empty vertices 
on the right.  Each empty vertex is connected (with a hedge) only to
the anti-vertex horizontally adjacent to it (see Figure
\ref{fig:graph}).

\begin{figure} [h]
\centering
\includegraphics{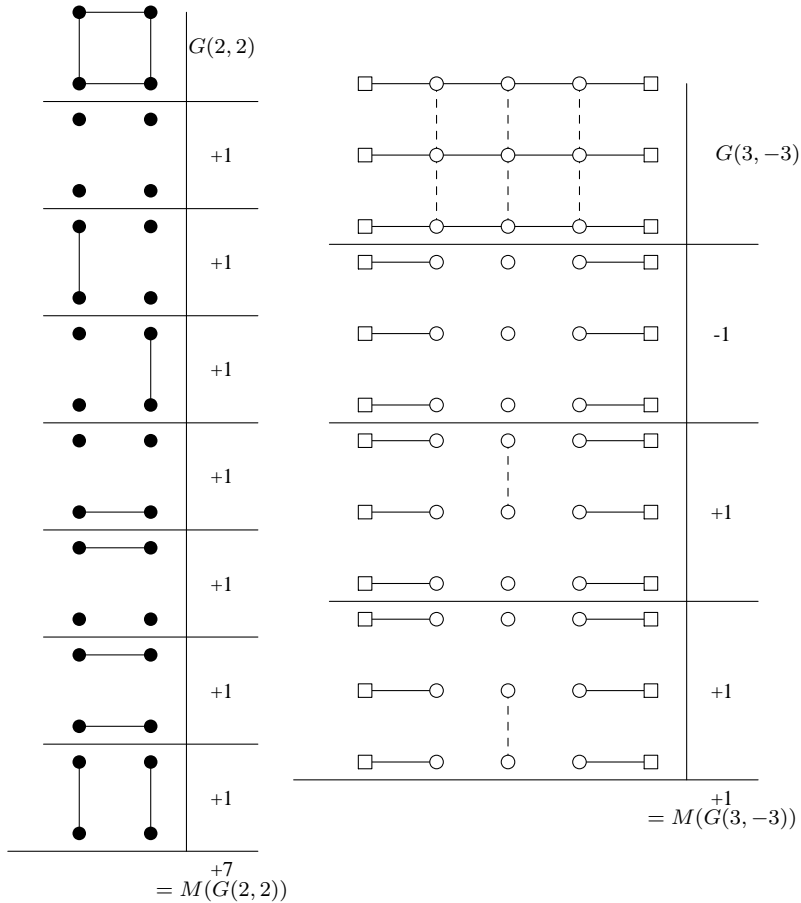}
\caption{Pictures of nonzero signed matchings of $G(2,2)$ and 
$G(3,-3)$.}
\label{fig:signed}
\end{figure}

For fixed $m > 0$ and any integer $n$, we define the number
$M(m,n) = M(G(m,n))$ to be the sum of the weights of the signed
matchings of $G(m,n)$.  In general, $M(G)$ is the sum of
the weights of all of the signed matchings of $G$. Some
examples are given in Figure \ref{fig:signed}.

Having made all our definitions in terms of rectangular grid
graphs, we would now like to extend them to apply to generalized
rectangles.  For any finite graph $G$, define the generalized 
rectangle graph $G\times P_{n}$ with $n > 0$ as follows:  
Picture $n$ copies of $G$ lined up side by side. Each vertex
in one copy of $G$ is connected with an edge to its image in the
copies of $G$ to the left and right of it. Define all
the edges \emph{within} a copy of $G$ to be vedges, and the edges
\emph{between} copies of $G$ to be hedges. With this idea, we can
define $G\times P_{-n}$ to be $G^{*}\times P_{n}$ with $m$ ($=$ 
the number of vertices in $G$) empty vertices to the left and to 
the right of it. Each empty vertex is connected to exactly one of 
the vertices of $G^{*}$ with a hedge.  We shall see that these
definitions are the right ones that will enable us to prove 
statements of reciprocity for more than just grid graphs.

\section{Adjunction}
\label{sec:adj}

Although it seems there are separate definitions for $G(m,n)$
(resp. $G\times P_{n}$) when $n > 0$ and when $n \leq 0$, we will
show they are actually the same object.  In Propp's words, they
fit together ``seamlessly." If they fit together this way, then it
follows that for any $G$, the sequence $M(G\times P_{n})$
satisfies a linear recurrence of finite degree with constant
coefficients. We will omit the proof as it is given in \cite{propp}.
The method of proof utilizes transfer matrix methods and the
Cayley-Hamilton theorem, and it requires only that the family of
graphs satisfy the property detailed below.

We begin by making a natural observation about the half-infinite
sequence $M(G(m,n))$, $n>0$, and then we will prove that it
actually holds for any $n$. In doing so, we will be able to
establish that we do indeed have an appropriate definition for
$G(m,n)$ (resp. $G\times P_{n})$ with $n \leq 0$, and all our
desires for the doubly infinite sequence---integrality,
uniqueness, and reciprocity---will be fulfilled.

\begin{figure} [h]
\centering
\includegraphics{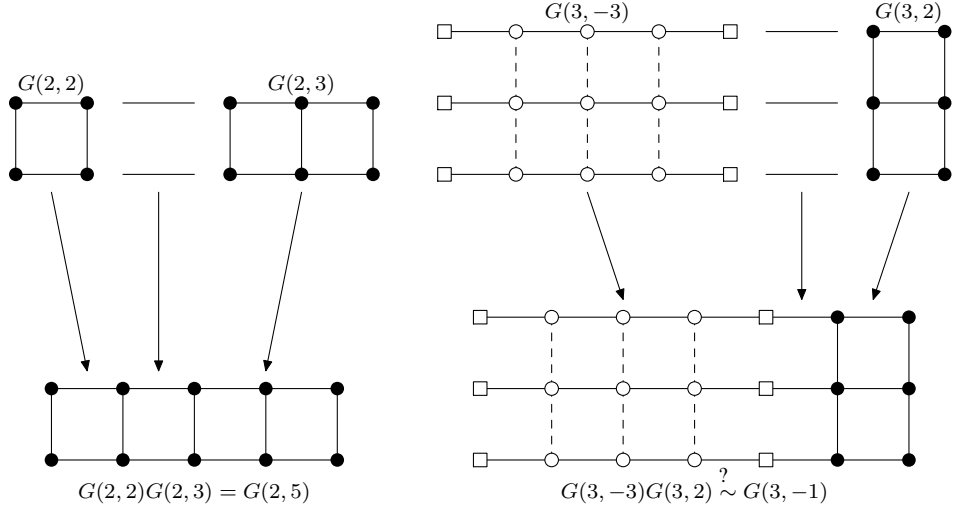}
\caption{Pictures of two pairs of adjoined graphs.}
\label{fig:adjoined}
\end{figure}

We first notice that for positive $n_{1}, n_{2}$, we can draw the
graph $G(m,n_{1}+n_{2})$ by placing $G(m,n_{1})$ and $G(m,n_{2})$
side by side and using hedges to connect the rightmost $m$
vertices of one to the leftmost $m$ vertices of the other, as
shown in Figure \ref{fig:adjoined}. Call this operation
\emph{adjunction}. Generally, define the adjunction of two graphs 
$H = G\times P_{n_{1}}$, $H'=G\times P_{n_{2}}$ to be a new graph 
formed by connecting $H$ to $H'$ with hedges, 
matching up corresponding vertices.  In particular:
\begin{itemize}
\item If $n_{1}, n_{2}$ are
positive, connect every vertex in the rightmost copy of $G$ in $H$ to 
its image
in the leftmost copy of
$G$ in
$H'$ using a hedge.
\item If $n_{1}$ is positive, $n_{2}$ negative, join every vertex
$v$ in the rightmost copy of $G$ in $H$ to an empty vertex on the left 
side of
$H'$ so that $v$ is connected to $v^{*}$ (the anti-vertex corresponding 
to $v$)
in the leftmost copy of $G^{*}$ in $H'$ by the path hedge-empty 
vertex-hedge.
(Similarly if $n_1$ is negative and $n_2$ is positive.)
\item If $n_{1}, n_{2}$ are both negative, join every vertex in the
rightmost copy of $G^{*}$ in $H$ to its image in the leftmost copy of 
$G^{*}$
in $H'$ by connecting their adjacent empty vertices with a hedge.
\end{itemize}
We write the adjunction of $G\times P_{n_{1}}$ and
$G\times P_{n_{2}}$ as $(G\times P_{n_{1}})(G\times P_{n_{2}})$.

Returning to case of rectangular grid graphs where $n_{1},n_{2}$ are 
positive,
$G(m,n_{1})G(m,n_{2}) = G(m,n_{1}+n_{2})$, and naturally,
$M(G(m,n_{1})G(m,n_{2})) = M(G(m,n_{1}+n_{2}))$. This is true more 
generally, 
as stated in the following:

\begin{theorem}[\bf{Adjunction}] \label{adj}
\[M((G\times P_{n_{1}})(G\times P_{n_{2}})\cdots (G\times P_{n_{k}})) =
M(G\times P_{n_{1} + n_{2} + \cdots + n_{k}})\] for all integers 
$n_{1},
\ldots, n_{k}.$  In particular, \[M(G(m,n_{1})G(m,n_{2})\cdots 
G(m,n_{k})) =
M(G(m, n_{1} + n_{2} + \cdots + n_{k}))\] for all integers $n_{1}, 
\ldots,
n_{k}$.
\end{theorem}

The proof of Theorem \ref{adj} will require the following two lemmas.

\begin{lemma}  Let $G$ be a signed graph where two vertices, $a$ and 
$b$, are
connected by the path hedge-empty vertex-hedge-empty vertex-hedge.
Then $G$ has same number of signed matchings as the graph $G'$, where
$G'$ is identical to $G$ except that $a$ and $b$ are connected
with one hedge.
\end{lemma}

\begin{figure} [h]
\centering
\includegraphics{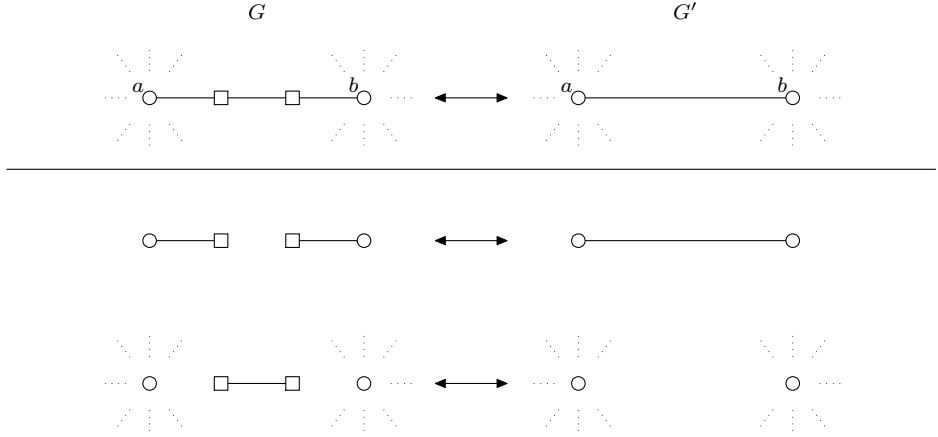}
\caption{Proof of bijection between $G$ and $G'$ in Lemma 1.}
\label{fig:lemma1}
\end{figure}

\emph{Proof:} The claim made by Lemma 1 is most easily seen
in Figure \ref{fig:lemma1} where we see the immediate bijection.
When a matching of $G$ does not contain the hedge joining the
two empty vertices, then the matching must contain the two
edges that join these two empty vertices to $a$ and $b$.
This corresponds to a matching of $G'$ that contains the hedge
joining $a$ and $b$.

On the other hand, any matching of $G$ where the hedge between the
empty vertices is present is a matching in which $a$ and $b$ are 
connected outwards.  Such a matching corresponds to a matching of
$G'$ that does not contain the hedge joining $a$ and $b$.
\qed

\begin{lemma}  Let $G$ be a finite graph with all plain components.  
Then let
$H$ be the graph $(G\times P_{n_1})(G\times P_{n_2})$ for $n_{1} >
0 > n_{2}$. Then $H$ has the same number of signed matchings as
the graph $H'$=$(G\times P_{n_1 -1})(G\times P_{n_2 +1})$.
\end{lemma}

\begin{figure} [h]
\centering
\includegraphics[scale=.8]{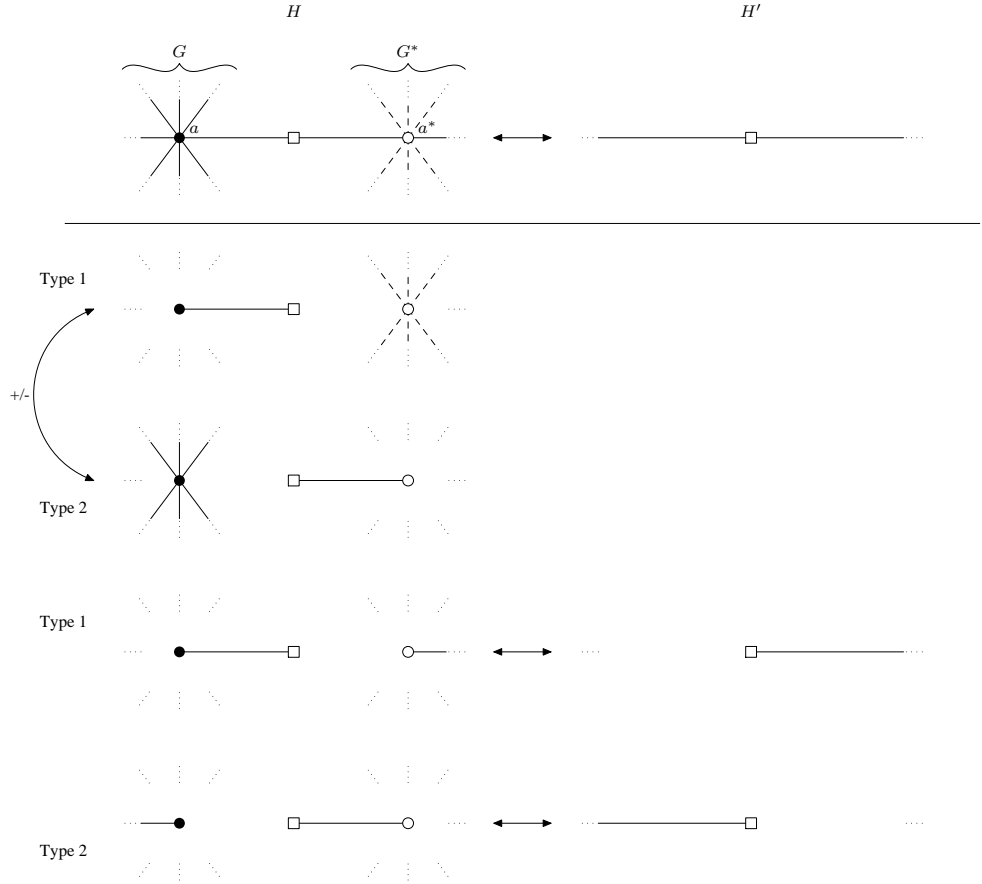}
\caption{Proof of Lemma 2.} \label{fig:lemma2}
\end{figure}

\emph{Proof:}
We begin by examining what happens at each vertex; see Figure
\ref{fig:lemma2}. Consider the plain vertex $a$ of $G$ connected
to the anti-vertex $a^{*}$ of $G^{*}$ by the path hedge-empty
vertex-hedge.  There are essentially two types of matchings of
$H$.  We will say type 1 matchings include the hedge connecting $a$
to the empty vertex and type 2 matchings include the hedge
connecting $a^{*}$ to the empty vertex.  Most of the matchings of
type 1 will cancel with most of the matchings of type 2.

Given a matching of type 1, there are only three cases for what
can happen the anti-vertex $a^{*}$ of $G^{*}$: it can be isolated,
it can have an anti-vedge incident with it, or it can have a hedge
incident with it.  Likewise for any matching of type 2, the vertex
$a$ can be isolated, incident with a vedge or incident with a
hedge. If there are $k$ signed matchings of type 1 where $a^{*}$
is isolated, then there are $-k$ signed matchings of type 2 where
$a$ is isolated.  Similar cancellation occurs between matchings of
type 1 where $a^{*}$ has an anti-vedge and matchings of type 2
where $a$ has a vedge (corresponding to the anti-vedge incident
with $a^{*}$).

The only remaining cases are those of type 1 and type 2 where $a$
and $a^{*}$ both have hedges.  We claim that these matchings are
in bijection with the graph where $a, a^{*}$, all their incident
vedges, and the hedges and empty vertex between $a, a^{*}$ are all
replaced with one empty vertex.  The correspondence is shown in
Figure \ref{fig:lemma2}.  Since $a$ was any vertex of $G$, this
can be done for every vertex of $G$ and the lemma holds.
\qed

With the lemmas proved, Theorem \ref{adj} is not difficult to
show.
\begin{figure} [h]
\centering
\includegraphics[scale=0.8]{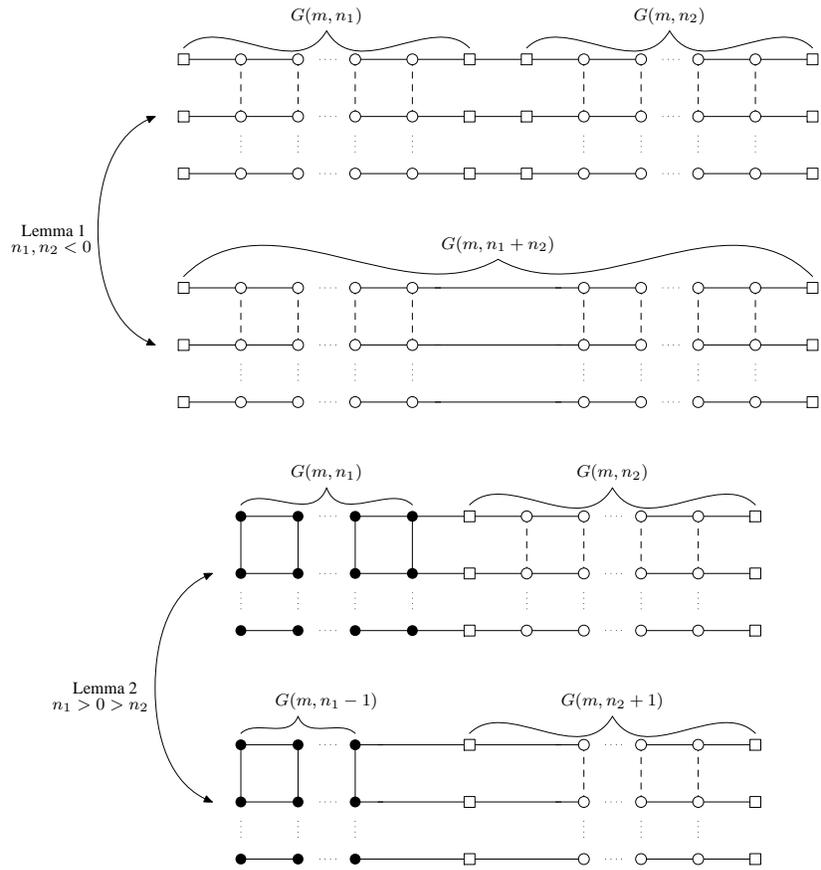}
\caption{Picture of Theorem \ref{adj} using lemmas.}
\label{fig:genrect}
\end{figure}

\emph{Proof of Theorem 1:} It is easily verified that $M((G\times
P_{0})(G\times P_n)) = M(G\times P_n)$, so we may assume all the
$n_{i}$ are non-zero. We have that $M((G\times P_{n_{1}})(G\times
P_{n_{2}})) = M(G\times P_{n_{1}+n_{2}})$ whenever $n_{1}$ and
$n_{2}$ are negative by applying Lemma 1 to each place where
adjoining takes place.  The case $n_1,n_2>0$ is even easier.
When $n_{1}$ and $n_{2}$ are of opposite sign, we can apply Lemma 2
repeatedly. Upon each application of Lemma 2 we change neither the
number of signed matchings nor the difference between the number of
vertices and the number of anti-vertices. \qed

\section{Reciprocity}
\label{sec:rec}

\subsection{Combinatorial Reciprocity}
\label{sec:crec}

The combinatorial statement of reciprocity is rather obvious given
the definition of $G\times P_{n}$ for $n \leq 0$ and the neediness
of empty vertices:

\begin{theorem}[\bf{Reciprocity (I)}] \label{recip1}
\[M(G\times P_{-n-2}) = M(G^{*}\times P_{n}).\] In particular,
\[M(G(m,-n-2)) = M(G^{*}(m,n)).\]
\end{theorem}

\begin{figure} [h]
\centering
\includegraphics[scale=1.25]{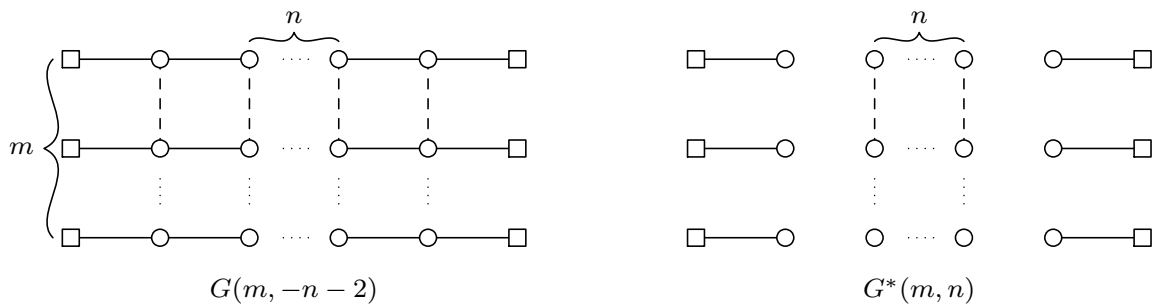}
\caption{Picture of reciprocity for grid graphs.}
\label{reciprocity}
\end{figure}

\emph{Proof:} On the left side of Figure \ref{reciprocity} we have
$G(m,-n-2)$. But since all the outside edges are forced and have
weight 1, the number of signed matchings of $G(m,-n-2)$ is clearly 
equal to the number of matchings of $G^{*}(m,n)$.  The situation 
for generalized rectangles is identical.
\qed

This statement of reciprocity is the natural extension of the
property of dimer coverings observed by Propp (\cite{propp}). Indeed,
if we were to consider only the perfect matchings of grid graphs,
then the statement would be the same: \[N(G(m,-n-2)) =
N(G^{*}(m,n)).\] It is known that any dimer covering can be
obtained from any other dimer covering by local moves that leave
the parity of vedges unchanged.  That is, the sign of every
\emph{perfect} matching of $G^{*}(m,n)$ is the same, so $|N(m,n)|
= |N^{*}(m,n)|$.  With less-than-perfect matchings, this is
clearly not the case.

Perhaps the strength of this result is more obvious when we give
all the edges and vertices formal weights and derive a statement
of reciprocity in terms of a generating function.  Though we could
do so for any finite graph of the form $G\times P_n$, we will only 
derive an explicit formula for $G(m,n)$.

\subsection{The Generating Function}
\label{sec:gfunc}

\begin{figure} [h]
\centering
\includegraphics[scale=.75]{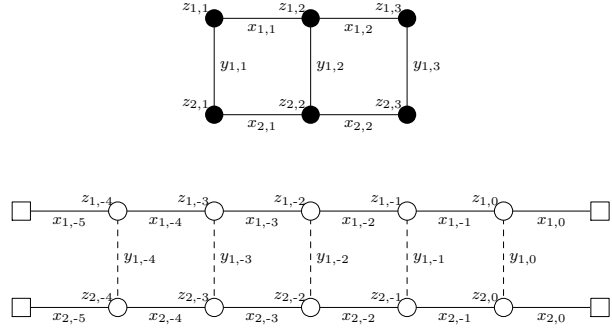}
\caption{Picture of graphs weighted with formal variables.}
\label{fig:weighted}
\end{figure}

For convenience, let the weights of the vertices of $G(m,n)$ for
$n>0$ be indexed by $\mathbb{N}\times\mathbb{N}$, read from left
to right and top to bottom.  Notice that the top left vertex is
$z_{1,1}$ and the bottom right vertex is $z_{m,n}$. Then $x_{i,j}$
is the weight given to the hedge between $z_{i,j}$ and
$z_{i,j+1}$, and $y_{i,j}$ is the weight given to the vedge
between $z_{i,j}$ and $z_{i+1,j}$. To get the weights for vertices
of $G(m,-n)$, we begin by giving the top right anti-vertex weight
$z_{1,0}$ and we proceed analogously. The weights of $G(2,3)$ and
$G(2,-5)$ are shown in Figure \ref{fig:weighted}.

The weight of a matching with formal variables is not
exactly like the weight as defined before, as is seen in the
following:
\begin{defin}
If $\mu$ is a matching of $G(m,n)$, define the weight of $\mu$,
$w(\mu)$, as follows (where $\mathrm{sgn}(t)=1$ if $t$ is positive,
$-1$ otherwise):
\[w(\mu) = \frac{\displaystyle\biggl(\prod_{\text{\rm hedges in }
\mu}\!\!\!\!\!x_{i,j}\biggr)\biggl(\prod_{\text{\rm vedges in }
\mu}\!\!\!\!\mathrm{sgn}(j)y_{i,j}\biggr)\biggl(\prod_{\text{\rm
vertices in }
\mu}\!\!\!\!\mathrm{sgn}(j)z_{i,j}\biggr)}{\displaystyle\prod_{\text{\rm
all hedges in } G(m,n)}\!\!\!\!x_{i,j}^{(1-\mathrm{sgn}(j))/2}}.\]
\end{defin}

\begin{figure} [h]
\centering
\includegraphics[scale=.75]{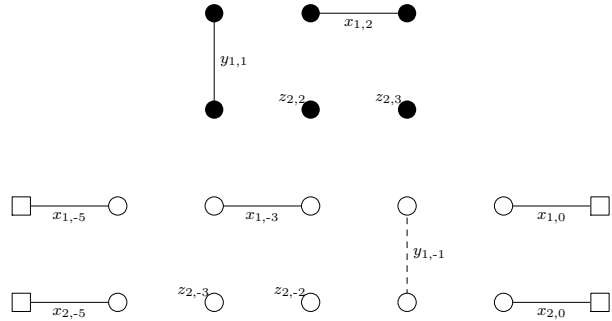}
\caption{Picture of matchings weighted with formal variables.}
\label{fig:weighted2}
\end{figure}

Notice that for $n<0$, there is a monomial denominator in the
weights of matchings. By our definition, the matchings shown in
Figure \ref{fig:weighted2} have weights
\[x_{1,2}y_{1,1}z_{2,2}z_{2,3} \mbox{ and }
\frac{y_{1,-1}z_{2,-2}z_{2,-3}}
{x_{1,-1}x_{1,-2}x_{1,-4}x_{2,-1}x_{2,-2}x_{2,-3}x_{2,-4}}\]

With this new definition of weights, we define the matching
polynomial. In the literature, there are conventions for
describing the matching polynomial of a graph. Our definition does
not adhere strictly to these conventions, though it is similar to
the partition function for the monomer-dimer model, with distinct
formal variables for all monomers and dimers (see
\cite{hl},\cite{ch}).

\begin{defin} The matching polynomial of $G(m,n)$,
$f_{n}(x_{i,j},y_{i,j},z_{i,j})$, is given by
\[f_{n}(x_{i,j},y_{i,j},z_{i,j}) 
= \sum_{\text{\rm all matchings $\mu$ of } 
G(m,n)}\!\!\!\!\!\!\!w(\mu)\]
where $i=1,\ldots,m$, $j=1,\ldots,n$ if $n>0$, $j=0,-1,\ldots,n$ if
$n<0$.
\end{defin}

Strictly speaking, the matching polynomial is a Laurent
polynomial: a polynomial in the variables $x_{i,j}$,
$x_{i,j}^{-1}$, $y_{i,j}$, $z_{i,j}$. Notice if we set all of the
weights equal to 1 we get $f_{n}(1,1,1) = M(G(m,n))$, the number
of signed matchings. Similarly if we set all the $x_{i,j},
y_{i,j}$ equal to 1 but set all the $z_{i,j}$ equal to 0, then
$f_{n}(1,1,0) = N(G(m,n))$, the number of signed perfect
matchings.

By construction the polynomials $f_{n}$ will satisfy a linear
recurrence very similar to that which governs $M(G(m,n))$. For
example, if $m =1$, then we get $f_{n} = z_{1,n}f_{n-1} +
x_{1,n-1}f_{n-2}$ for all $n$.  Likewise, there is a link between
the polynomials $f_{n}$ and $f_{-n-2}$ for $n\geq 0$, seen most
easily if we let $x_{i,j}=x$, $y_{i,j}=y$, $z_{i,j}=z$ for all
$i,j$:
\begin{equation}
f_{n}(x,-y,-z) = x^{m(n+1)}f_{-n-2}(x,y,z).
\end{equation}

As mentioned in Section 3, for any fixed $m$, the numbers
$M(m,n)$, $n \in \mathbb{Z}$ satisfy a linear recurrence of finite
degree with constant coefficients.  Therefore there is a rational
generating function for the number of weighted matchings,
\[F_{m}(t,x,y,z) = \sum_{n\geq 0}f_{n}(x,y,z)t^{n}.\] Given (1) it
is an easy exercise (see for example \cite{stan1}, ch. 4) to state
a reciprocity theorem for the generating function:

\begin{theorem}[\bf{Reciprocity(II)}]\[x^{m}t^{2}F_{m}(t,x,y,z) =
-F_{m}\biggl(\frac{1}{tx^{m}},x,-y,-z\biggr).\]
\end{theorem}

This completes the main goal of the paper.

\section{More on Linear Recurrences and Reciprocity}
\label{sec:more}

Though we have only proved theorems here for graphs of the type
$G\times P_{n}$ (cylinders of fixed circumference and varying
height for example), we have also been able to use a modified form
of adjunction to build a model for the honeycomb graph of fixed
height and width $n \in \mathbb{Z}$. Other graphs that we examined
but were unable to apply our methods to include cylinders with fixed
height and varying circumference, M\"obius strips, tori, and
projective planes. We do however feel that there may be a way of
dealing with such graphs.

In his unpublished paper \cite{speyer}, David Speyer developed a 
matrix method for encoding perfect matchings of a graph.  With this
method he was able to state theorems about recurrences and
reciprocity for a broader range of graphs than those we have
handled here.  In particular, he was able to make a statement
about the M\"obius strip and projective plane. An adaptation of
his method to general matchings seems promising, though
experimentation has shown it to be less than straightforward.

\section{Acknowledgements}
\label{sec:ack}

We'd like to thank Jim Propp for introducing us to this problem
and for getting us started with the early investigations.  Without
his REACH program to bring the authors together, this work could
not have been done.  Thanks also to Harvard University, Brandeis
University, the University of Massachusetts-Boston, the University
of Wisconsin, the NSF, and the NSA.

\end{document}